\documentclass[12pt,a4paper]{article}
\usepackage[utf8]{inputenc}
\usepackage[english,russian]{babel}
\usepackage{amsmath}
\usepackage{amssymb}
\usepackage{amsfonts}
\usepackage{amsthm}
\usepackage[matrix,arrow]{xy}
\usepackage{longtable}
\usepackage{graphicx}

\newtheorem{theorem}{Theorem}[section]
\newtheorem{stm}[theorem]{Statement}

\newtheorem{lemma}[theorem]{Lemma}
\newtheorem{imp}[theorem]{Corollary}
\newtheorem{conj}[theorem]{Conjecture}

\theoremstyle{definition}
\newtheorem{df}[theorem]{Definition}

\newcommand{\mb}{\mathbb}
\newcommand{\ra}{\rightarrow}
\renewcommand{\l}{\lambda}

\renewcommand{\proof}{{\it Proof.}\;\;}

\begin{document}

\begin{center}
{\Large\textbf{On the number of coverings of the sphere\\ ramified over given points}}\\
{\textbf{Boris~Bychkov}}\footnote{Research is supported by the RFBR grants 12-01-31233, 13-01-00383}\\
{\small\emph{Department of mathematics, National Research University Higher School of Economics, Vavilova str. 7, 117312, Moscow, Russia }}\\
{\small\emph{bbychkov@hse.ru}}
\end{center}

\noindent\textbf{Abstract}

We present the generating function for the numbers of isomorphism classes of coverings of the two-dimensional sphere by the genus $g$ compact oriented surface not ramified outside of a given set of $m+1$ points in the target, fixed ramification type over one point,
and arbitrary ramification types over the remaining $m$ points. We present the genus expansion of this generating function and prove, that the generating function of coverings of genus $0$ satisfies some system of differential equations. We show that this generating function is a specialization of the function from paper \cite{GJ} and, therefore, satisfies the KP-hierarchy.

\section{Introduction}

The problem of enumerating coverings of the sphere by two-dimensional surfaces with fixed
ramification types over given points was posed by Hurwitz~\cite{Hu}.
During the last decades, we saw growing interest to this problem, due to
discovery of its connections with various physical theories, Gromov--Witten invariants and geometry of moduli spaces of complex curves. In this paper, we consider a close problem, namely, that of enumerating isomorphism classes
of coverings of the two-dimensional sphere by the genus $g$ compact oriented surface
not ramified outside of a given set of $m+1$ points in the target, fixed ramification type over one point,
and arbitrary ramification types over the remaining $m$ points. We start with paper \cite{BMS} by M.~Bousquet-Melou and G.~Schaeffer, where the following formula for the number of genus~$g=0$ such coverings was deduced.

{\theorem[Bousquet-Melou, Schaeffer]  Let $\sigma_0\in S_n$ be a permutation having $d_i$ cycles of length $i\;(i=1,2,3,\dots)$ and let $l(\sigma_0)$ be the total number of cycles in $\sigma_0$, then
\begin{equation}
G_{\sigma_0}(m) = m\dfrac{(mn-n-1)!}{(mn-n-l(\sigma_0)+2)!}\prod\limits_{i\geq 1}\left(i\binom{mi-1}{i}\right)^{d_i}.\label{0}
\end{equation}}

In this paper, we present the generating function for the numbers of isomorphism
classes of coverings of arbitrary genus. Our argument is based on M.~Kazaryan's remark about the
operator acting in the center of the group algebra of a symmetric group
by multiplication by the sum of all permutations and its eigenvalues.

\section{Results}
\hspace{6mm}Let $S_n$ be the group of permutations of $n$ elements. Let $\nu_1,\ldots,\nu_t$ be the
lengths of the cycles in the decomposition of a permutation $\sigma \in S_n$ into the product of disjoint cycles; we say that $\sigma$ {\it has the cycle type $\nu$}. Introduce notation $l(\nu) = t,\; n = |\nu| = \nu_1 + \cdots + \nu_t$ and let $\nu = \rho(\sigma) = \{\nu_1,\ldots,\nu_t\}$ be the set of the lengths of the cycles.

Irreducible representations of $S_n$ are in one-to-one correspondence with partitions $\nu\vdash n$ of $n$
(see, e.g., \cite[\S 3]{Na}). We will denote irreducible representa-\\tions by partitions.

Denote by $\mathrm{dim}_{\nu}$ the dimension of the irreducible representation $\nu$ of $S_{|\nu|}$.

Denote by $|\mathrm{Aut}(\nu)|$ the order of the automorphism group of the partition $\nu$.
If the partition $\nu$ has $d_i$ parts of length $i,\; (i=1,2,3,\ldots)$, then
$|\mathrm{Aut}(\nu)| = d_1!d_2!d_3!\ldots$.

Denote by $p_i = p_i(x_1,\ldots,x_n)$ the symmetric power polynomial in $n$ variables defined by the formula
\begin{equation}
p_i(x_1,x_2,\ldots,x_n) = \sum\limits_{j=1}^n x_j^i.\label{7}
\end{equation}

Denote by $s_\nu$ the {\it Schur function} corresponding to the partition $\nu$.
Schur functions are polynomials in the symmetric power polynomials $p_i$ (Definition~\ref{d1} below).
Define the corresponding {\it scaled Schur function} by the equation
\begin{equation*}
s\hbar_\nu(\hbar,p_1,p_2,p_3,\ldots) = s_\nu(p_1,p_2\hbar,p_3\hbar^2,\ldots).
\end{equation*}

For a partition $\nu = \{\nu_1,\ldots,\nu_t\}$, denote by $p_\nu$ the product $p_\nu = p_{\nu_1}\ldots p_{\nu_t}$.

\vspace{4mm}
Recall that partitions can be usefully presented by intuitive geometric objects, {\it Young diagrams}. A Young diagram is a finite subset of the two-dimensional positive integer quadrant such that together with any point
of the integer lattice it contains all the lexicographically smaller points.
Points are usually represented by unit squares, the total amount of points is equal to $|\nu|$,
the squares are arranged in downright rows of nonincreasing lengths $\nu_1\ge\nu_2\ge\dots\ge\nu_t$,
 and the total number~$t$ of the rows is equal to the size $l(\nu)$ of the partition $\nu$.

\df The {\it content} of a cell $k$, which lies on the intersection of the $i$th column and $j$th row of the
Young diagram, is the value $c(k) = j - i$. For an arbitrary cell in the Young diagram,
 consider the set of cells consisting of itself and all the cells that lie
 in the same row to the right of it and in the same column down of it. This set is the {\it hook} of the cell~$k$.
  The {\it length of the hook} $h(k)$ is the number of cells in this hook.

\vspace{4mm}
Let us recall the connection between ramified coverings of the two-dimensio\-nal sphere and decompositions of a permutation into products of permutations.

\df Let $X$ and $Y$ be two topological spaces, $Y$ being path connected, let $f:X\rightarrow Y$
 be a continuous mapping, and let~$S$ be a discrete set.  If for any point $y\in Y$
 there exists a neighbourhood $V$ of $y$ such that the preimage $f^{-1}(V)\subset X$
  is homeomorphic to $V\times S$,
  then the triple $(X,Y,f)$ is called an \textit{unramified covering} of $Y$ by $X$
  with the fiber~$S$. The cardinality of the set $S$ is called
   the {\it degree} of the covering. A covering of degree~$n$ is also said to be {\it $n$-sheeted}.

Let $f:X\rightarrow Y$ be an unramified covering, $y_0\in Y$. Every continuous path $\gamma : [0,1] \rightarrow Y$ that starts and ends in $y_0$ defines a permutation of the set $f^{-1}(y_0)$. This permutation is called the {\it  monodromy\/}
along~$\gamma$.

Let $f:X\rightarrow\mathbb{C}P^1\setminus T$ be a finite-sheeted covering of
$\mathbb{C}P^1$ with $k$ punctures $y_1,\ldots, y_k,\; T = \{y_1,\ldots, y_k\}\subset\mathbb{C}P^1$.
The projective line $\mathbb{C}P^1$ is endowed with the orientation induced by the complex structure.
Let $y_0 \in \mathbb{C}P^1\setminus T$. Consider $k$ oriented paths $c_i$ connecting $y_0$ with $y_i,\;i = 1,\ldots, k$ on $\mathbb{C}P^1$ such that they do not have intersections outside $y_0$ and enter $y_0$ in the order prescribed
by their numbering. Transform every path $c_i$ into a loop $\gamma_i \in \pi_1(\mathbb{C}P^1\setminus T,y_0)$
going around the point~$y_i$ in the positive direction an returning back along~$c_i$ to~$y_0$.
We denote by~$g_i$ the monodromy permutation along the loop $\gamma_i$. The group $G=\langle g_1,\ldots,g_k\rangle$ generated by monodromy permutations is called \textit{the monodromy group\/} of the covering.
It consists of all the monodromy permutations along paths starting and ending in~$y_0$.

To each finite-sheeted covering of a punctured sphere, the unique ramified covering of the sphere without punctures
is associated: compactify the punc-\\tured sphere $\mathbb{C}P^1\setminus T$ by adding all the points $y_1,\ldots, y_k$ and for every point $y_i$ add to $X$ as many points as there are independent cycles in the permutation $g_i$. New points in $X$ will be the preimages of the new points $y_i$ with multiplicities equal to the lengths of the corresponding cycles in the permutation $g_i$.

\df Let $X$ be a compact two-dimensional orientable surface. If there exists a finite set of points $T = \{y_1,\ldots, y_k\}\subset \mathbb{C}P^1$ such that $f$ is obtained from an unramified covering of the punctured sphere $\mathbb{C}P^1\setminus T$ by the construction described above, then the continuous mapping
 $f:X\rightarrow\mathbb{C}P^1$ is called a \textit{ramified covering}. The {\it ramification type\/} of a ramification point $y_i$ is the cycle type of the permutation $g_i$.
\\
Let $f_1:X_1\ra \mb{C}P^1$ and $f_2:X_2\ra \mb{C}P^1$ be two ramified coverings. If there exists an orientation
preserving homeomorphism $u:X_1\ra X_2$ such that the diagram below is commutative,
then the two ramified coverings are said to be {\it isomorphic}:
$$
\xymatrix{
X_1 \ar[dr]_{f_1} \ar[rr]^u && X_2 \ar[dl]^{f_2}\\
& \mb{C}P^1
 }
$$

\vspace{4mm}
Denote by $b_{g,\nu,m}$ the total number of isomorphism classes of
ramified coverings of the two-dimensional sphere by a surface of genus $g$ having $m$ ramification points
with arbitrary ramification types and a distinguished ramification point with given ramification type $\nu = \{\nu_1,\ldots,\nu_c\}$ over it. Denote by $S$ the following generating function for the numbers $b_{g,\nu,m}$:
\begin{equation*}
S(\hbar,p_1,p_2,\ldots;m) = \sum\limits_{n=0}^{\infty}\sum\limits_{\nu\vdash n} b_{g,\nu,m}p_\nu \hbar^{2g},
\end{equation*}
where $p_\nu = p_{\nu_1}p_{\nu_2}\ldots$ and the genus $g$ of the covering surface
   can be computed from the Riemann--Hurwitz formula, $2-2g = 2n - \sum\limits_P (k(P)-1)$,
   $k(P)$ is the ramification order at the ramification point~$P$, see Sec.~\ref{s4}.

\vspace{4mm}
The main result of this paper is the following
{\theorem The function $S$ admits the following representation:
\begin{equation}
S(\hbar,p_1,p_2,\ldots;m) = \hbar^2\log \left(\sum\limits_{n=0}^{\infty}\sum_{\nu\vdash n} \prod\limits_{k\in\nu}(1+c(k)\hbar)^m \dfrac{\mathrm{dim}_\nu}{n!} \hbar^{-2n}s\hbar_\nu \right). \label{1}
\end{equation}\label{t1}
}

The generating function $S$ has the genus expansion
\begin{multline*}
S(\hbar,p_1,p_2,\ldots;m) = S_0(p_1,p_2,\ldots;m) + \hbar^2S_1(p_1,p_2,\ldots;m) + \\
\hbar^4S_2(p_1,p_2,\ldots;m)+\ldots,
\end{multline*}
where the functions $S_g,\; g=0,1,2\ldots$, are the generating functions
for numbers of coverings of genus g. The coefficients $b_{0,\nu,m}$ of the function $S_0$ are the Bousquet-Melou--Schaeffer numbers $G_{\sigma_0}(m)$ (see Eq.~(\ref{0})) divided by $|\mathrm{Aut}(\nu)|\prod\limits_{i=1}^{l(\nu)}\nu_i$,
 where $\nu$ is the cycle type of the permutation $\sigma_0$. Indeed, the Bousquet-Melou--Schaeffer number is equal to the number of genus~$0$ decompositions of a given
permutation $\sigma_0$, while $b_{0,\nu,m}$ enumerate decomposi-\\tions of arbitrary permutations from a given conjugacy class.

For higher genera, Eq.~(\ref{1}) produces an efficient way to compute the coefficients of the expansion. For example,
the generating function $S_1$ for genus~1 (the coefficient of the series (\ref{1}) of $\hbar^2$) up to the covering degree~$4$ is

\begin{multline*}
S_1(p_1,p_2,\ldots;m) = \dfrac{1}{48}m(m-1)(m-2)(m-3)p_1^2 + \dfrac{1}{12}m(m-1)(m-2)p_2 + \\
\dfrac{1}{72}m(m-1)(m-2)(4m^3-21m^2+35m-20)p_1^3 + \\
\dfrac{1}{6}m(2m-3)(m-2)(m-1)^2p_2p_1 + \dfrac{1}{24}m(3m-5)(m-1)(3m-2)p_3 + \\
\dfrac{1}{96}m(m-1)(m-2)(13m^5-99m^4+297m^3-445m^2+337m-105)p_1^4 + \\
\dfrac{1}{24}m(m-1)^2(m-2)(26m^3-103m^2+135m-60)p_2p_1^2 + \\
\dfrac{1}{24}m(m-1)^2(4m-5)(4m^2-10m+5)p_2^2 + \\
\dfrac{1}{16}m(m-1)^2(3m-5)(3m-4)(3m-2)p_3p_1 + \\
\dfrac{1}{12}m(m-1)(4m-3)(2m-1)(2m-3)p_4+\dots.
\end{multline*}

\vspace{4mm}
Generating series enumerating coverings often are solutions to integrable hierarchies, see, e.g.,~\cite{Ok}. The generating function $S$ is a specialization of a function from paper \cite{GJ} and, therefore, is a solution to the KP-hierarchy. 
In particular, it satisfies the first of the infinite series of equations in the KP-hierarchy:
\begin{equation*}
\dfrac{\partial^2 S}{\partial p_1\partial p_{3}} - \dfrac{\partial^2 S}{\partial p_2^2} = \dfrac{1}{2}\dfrac{\partial^2 S}{\partial p_1^2} - \dfrac{1}{12}\dfrac{\partial^4 S}{\partial p_1^4}.
\end{equation*}
Let $Y(\nu) = \prod\limits_{k\in\nu}y_{c(k)}$ denote the content product for the Young diagram corresponding to the partition $\nu\vdash n$, for the indeterminates $y_c$, $c=\dots,-2,-1,$ $0,1,2,\dots$.
Denote by $F$ the generating function
\begin{equation*}
F(\dots,y_{-2},y_{-1},y_{0},y_{1},y_{2},\dots;p_1,p_2,\ldots)=\log\left(\sum_{n=0}^\infty\sum_{\nu\vdash n}\prod_{k\in\nu}y_{c(k)} \frac{dim_\nu}{n!}p_\nu\right).
\end{equation*}
The following statement is the immediate consequence of Theorems 2.3 and 3.1 from \cite{GJ}:

{\stm The generating function $F$ is a solution to the KP-hierarchy.}

Consider the perturbation of the KP-hierarchy by the substitution $p_i = \dfrac{p_i}{\hbar^{i+1}}$,
where $\hbar$ is a formal parameter and $i=0,1,2,\dots$. For example the first equation of the hierarchy will be
\begin{equation*}
\dfrac{\partial^2 S}{\partial p_1\partial p_{3}} - \dfrac{\partial^2 S}{\partial p_2^2} = \dfrac{1}{2}\dfrac{\partial^2 S}{\partial p_1^2} - \dfrac{\hbar^2}{12}\dfrac{\partial^4 S}{\partial p_1^4}.
\end{equation*}

{\imp The function $S(\hbar,p_1,p_2,\ldots;m)$ is a solution to the perturbed KP-hierarchy.}

\proof
One can obtain the function $S$ from $F$ by the substitutions $y_{c} = (1+c\hbar)^m$ and $p_i = \dfrac{p_i}{\hbar^{i+1}}$.
\qed

\vspace{4mm}
Computations using Theorem~\ref{t1} lead to the following conjectural version of the Bousquet-Melou--Schaeffer
 formula to the case of genus $g=1$.

{\conj The Bousquet-Melou--Schaeffer numbers for coverings of the sphere by the torus admit the following representation:
\begin{equation}
b_{1,\nu,m} = P_{2t-1}(m,\nu)m\prod\limits_{i=1}^t (m\nu_i-2)_{(\nu_i-1)}.
\end{equation}
 Here $P_{2t-1}$ is a polynomial of degree $2t-1$ and $(m\nu_i-2)_{\nu_i-1} = (m\nu_i-2)(m\nu_i-3)\ldots (m\nu_i-\nu_i)$
 is the descending factorial.
}

In Sec.~3, we present necessary well-known facts about symmetric functions and representations of symmetric groups.
In Sec.~4, we prove main statements of the paper.

\section{Representations of symmetric groups}
\hspace{6mm}Let us reformulate our problem in the language of permutations. We will follow the terminology of~\cite{LZ}.

\df A {\it constellation} is a sequences of permutations $[g_1,\ldots, g_k],$ $g_i \in S_n$ such that the group $\langle g_1,\ldots, g_k \rangle$ acts transitively on an~$n$-element set and $g_1\cdot\ldots\cdot g_k = id$.
It is not hard to check that the $k$-tuple of permutations $[g_1,\ldots, g_k]$
generating the monodromy group of an unramified covering of a punctured sphere form a constellation.
This construction works in the opposite direction as well: for any constellation,
there exists a corresponding unramified covering of the punctured sphere.

Consider two constellations $[g_1,\ldots, g_k]$ and $[g'_1,\ldots, g'_k]$.
If there exists a permutation $h \in S_n$ such that $g'_i = h^{-1}g_ih$ for all $i = 1,\ldots, k,$ then these two constellations are said to be {\it isomorphic}. Two unramified coverings of a punctured
sphere $\mathbb{C}P^1\setminus Y$ are isomorphic if and only if the corresponding constellations are isomorphic.
{\theorem[\rm\bf Riemann's existence theorem] Consider any sequence of points $[y_1,\ldots, y_k]\in \mathbb{C}P^1$ and any constellation $[g_1,\ldots, g_k],\; g_i \in S_n$. Then there exists a Riemann surface~$X$ and a meromorphic function $f:X\rightarrow \mathbb{C}P^1$ such that $y_1,\ldots, y_k$ are the points of ramification of~$f$ and $g_1,\ldots, g_k$
are the corresponding monodromy permutations. The unramified covering $f:X\rightarrow \mathbb{C}P^1$ is unique up to isomorphism.}

One can find a proof in \cite{LZ}.

Thus, the numbers $b_{g,\nu,m}$ enumerate decompositions of permutations $\sigma_0\in S_n$ with given cycle type $\nu$ into a product of~$m$ permutations from $S_n$ (we count $m$-tuples of permutations up to common conjugation) provided the group generated by these permutations acts transitively on the set of $n$ elements. The genus of the covering surface
can be determined from the Riemann--Hurwitz formula.

\df The {\it group algebra} $\mathbb{K}G$ of a finite group $G$ is the
 $|G|$-dimensional vector space over the field $\mathbb{K}$ freely spanned by the elements of $G$.
 The product in $\mathbb{K}G$ is induced by the group operation in~$G$.

We will use only the field $\mathbb{K} = \mathbb{C}$.

Proofs of all
the facts below about the group algebras of symmetric groups
can be found, e.g., in~\cite{Vi} and~\cite{Mc}.

It is well known that every linear representation~$R:G\to GL(V)$ of a group~$G$
in a vector space $V$ over~$\mathbb{C}$ admits the unique extension to a linear representation
of the algebra $\mathbb{C}G$ in the same space according to the formula
$R(\sum\limits_{g\in G}a_gg) = \sum\limits_{g\in G}a_gR(g)$. Note that the inverse statement
is also true, so that there is a natural bijection between representations of a group $G$
and those of the group algebra $\mathbb{C}G$.

Let~$T$ be the regular representation of the algebra $\mathbb{C}G$,
 that is, the representation of $\mathbb{C}G$ on itself (if we think about $\mathbb{C}G$ as about a vector space)
  defined by the rule $T(a)x = ax,\; \forall a,x\in \mathbb{C}G$.

Define the scalar product in $\mathbb{C}G$ in the standard way:
\begin{equation}
(a,b) = \mathrm{tr}\,T(ab) = \mathrm{tr}\,T(a)T(b).\label{d}
\end{equation}

Denote by $C_\nu$ the element in the group algebra~$\mathbb{C}S_{|\nu|}$ of the symmetric group equal to the sum of all permutations whose cycle type is~$\nu$. The elements $C_\nu$ form
 a basis of the center of the group algebra $\mathbb{C}S_{|\nu|}$.

Consider the space $\mathbb{C}[G]$ of functions on a group~$G$. The formula
\begin{equation}
\varphi(\sum\limits_{g\in G}a_gg) = \sum\limits_{g\in G}a_g\varphi(g),\; \varphi\in\mathbb{C}[G], \label{5}
\end{equation}
extends every function $\varphi\in\mathbb{C}[G]$ to a linear function on the group algebra $\mathbb{C}G$.
Hence we have a natural bijection between $\mathbb{C}[G]$ and the dual space $\mathbb{C}G^*$
to the space $\mathbb{C}G$ given by Eq.~(\ref{5}).

The scalar product~(\ref{d}) defines an isomorphism between $\mathbb{C}G$
and its dual space, $g\mapsto\varphi_g,$ where
\begin{equation}
\varphi_g(h) = (g,h) =
\left\{\begin{aligned}
|G|,\; gh = e\\
0,\; gh \neq e\\
\end{aligned}
\right.
\;\;\forall g,h\in\mathbb{C}G.\label{d2}
\end{equation}

Let us transfer the scalar product to the space $\mathbb{C}[G]$:
\begin{equation}
(\varphi_1,\varphi_2) = \dfrac{1}{|G|}\sum\limits_{g\in G}\varphi_1(g)\varphi_2(g^{-1}).\label{d3}
\end{equation}

Let $R:G\ra GL(V)$ be an arbitrary linear representation of a group $G$.
\df Define the function $\chi\in \mb{C}[G]$ by the formula $\chi(g) = \mathrm{tr}\,R(g),\; g\in G$.
This function is called the \textit{character} of the representation $R$.

Denote by $\chi^\nu$ the character of the irreducible representation of $S_n$ \\
corresponding to the partition $\nu\vdash n$. Characters $\chi^\nu$ are idempotents:
\begin{equation}
\chi^\mu\chi^\nu = \dfrac{\mathrm{dim_\nu}}{n!}\delta_{\mu}^{\nu}\chi^\mu.
\end{equation}

Consider the mapping $\psi$ from $S_n$ to the space of quasihomogeneous polynomials of degree $n$ in variables  $p_i$, $\psi:\;\sigma\mapsto p_{\rho(\sigma)} = p_{\nu_1}\ldots p_{\nu_t}$, where weight $i$ is assigned to the variable $p_i$.

Define the \textit{characteristic mapping} $\mathrm{ch}$ from $Z\mb{C}S_n^*$ to the space of quasiho-\\mogeneous polynomials of degree~$n$ in the variables~$p_i$ by the formula:
\begin{equation*}
\mathrm{ch}(f) = \dfrac{1}{n!}\sum\limits_{g\in G}f(g)\psi(g).
\end{equation*}

\df Let $\nu$ be a partition of length less than $l+1$, then the
\textit{Schur function} $s_\nu(x_1,x_2,\ldots,x_l)$ of $\nu$ is the quotient of two determinants:\label{d1}
\begin{equation}
s_\nu(x_1,x_2,\ldots,x_l) = \dfrac{\mathrm{det}(x_i^{\nu_j+l-j})_{1\leq i,j\leq l}}{\mathrm{det}(x_i^{l-j})_{1\leq i<j\leq l}}.
\end{equation}

We will use Schur functions rewritten in the variables $p_1,p_2,\ldots$.
Their equivalent definition looks like follows.
First, the Schur function $s_k(p_1,p_2,\ldots)$ of a one-part partition is the coefficient of $t^k$ in the series $\mathrm{exp}\left(\sum\limits_{i=1}^\infty\dfrac{p_i}{i}t^i\right)=\sum\limits_{i=0}^\infty s_kt^k$.
Next, the Schur function $s_\nu(p_1,p_2,\ldots)$ is the determinant of the following matrix formed by one-part Schur functions:
$$
s_\nu(p_1,p_2,\ldots) = \mathrm{det}(s_{\nu_i-i+j})_{1\leq i,j\leq l(\nu)}.
$$

{\stm[\cite{Mc},1.7.3] The mapping $\mathrm{ch}$ is an isometric bijection between
the center of the group algebra $Z\mb{C}S_n^*$ and the space of quasihomogeneous polynomials of degree $n$,
under the isomorphism $\mathrm{ch}(\chi^\nu) = s_\nu$.}

In addition, there is a correspondence between $Z\mb{C}S_n$ and $Z\mb{C}S_n^*$.
Thus, the element $C_\nu=C_{\nu_1,\ldots,\nu_t} \in Z\mb{C}S_n$
is taken to the monomial $|C_\nu|p_{\nu_1}\ldots p_{\nu_t}$,
where $|C_\nu|$ is the number of elements in the conjugacy class
of a permutation of cycle type $\nu$.

\section{Proof of the main theorem}\label{s4}
\hspace{6.2mm}In this section we will prove the main Theorem \ref{t1}.

Any element $a\in Z\mb{C}S_n$ has expansion in the both bases $C_\nu$ and $\chi^\nu$:
$$
a = \sum\limits_{\nu\vdash n} \dfrac{(a,C_\nu)}{|C_\nu|n!}C_\nu = \sum\limits_{\nu\vdash n} (a,\chi^\nu)\chi^\nu.
$$
Due to the fact that the characters  $\chi^\nu$ are idempotent, we have
$$
a\cdot b = \sum\limits_{\nu\vdash n} (a,\chi^\nu)(b,\chi^\nu)\chi^\nu,\; \forall\; a, b\in Z\mb{C}S_n.
$$

Let us associate to an arbitrary element $a\in Z\mb{C}S_n$
the operator on $Z\mb{C}S_n$ acting by multiplication by $a$.
Then $\chi^\nu$ are eigenvectors of this operator with the eigenvalues $(a,\chi^\nu).$

Consider the operator $B:Z\mb{C}S_n\ra Z\mb{C}S_n$ defined by the formula $B = \sum\limits_{\nu \vdash n} \hbar^{|\nu|-l(\nu)} C_\nu$ (in other words, $B$ is a scaled sum of all the elements of the group $S_n$); here $\hbar$ is a formal variable.
Recall that we count $m$-tuples of permutations whose product is
a permutation of cyclic type $\nu$. Expand the operator $B^m$ in the eigenbasis of characters and
note that the eigenvalue of the vector $\chi^\nu$ is exactly the desired number.
Denote by $B_\nu$ the eigenvalue of the operator $B$ on the eigenvector $\chi^\nu$.
When computing~$B_\nu$, we will use the following statement.

{\lemma[\cite{KOV},4.1.2] The function $\Gamma\in \mb{C}[S_n], \Gamma:\sigma_\nu\mapsto \hbar^{l(\nu)}$, where $\hbar$ is a formal variable and $\sigma_\nu$ is an arbitrary permutation of cyclic type $\nu$,
has the following expansion in the basis of characters:
\begin{equation*}
\Gamma(\cdot) = \sum\limits_{\nu\vdash n} \prod\limits_{k\in\nu}\dfrac{\hbar+c(k)}{h(k)}\chi^\nu;
\end{equation*}
here $k$ runs over the set of cells of the Young diagram corresponding to the partition $\nu$,
$c(k)$ is the content of the cell $k$, $h(k)$ is the length of the corresponding hook.\label{l}
}

For the sake of completeness, we reproduce the proof of this statement from~\cite{KOV}.

\proof
By definition (see \cite{Mc}, 1.7), for an arbitrary function $f\in \mb{C}[ZS_n]$, the equality $\mathrm{ch}(f) = \sum\limits_{\nu\vdash n}z_\nu^{-1}f(\sigma_\nu)p_\nu$ is valid, where $z_\nu = 1^{\nu_1}2^{\nu_2}\ldots\nu_1!\nu_2!\ldots$.

Consider the generating series $\sum\limits_{n\geq 1} \mathrm{ch}(\Gamma)u^n$ for the numbers $\mathrm{ch}(\Gamma)$. By definition of the functions $\mathrm{ch}$ and $z_\nu$, we have
\begin{equation*}
1 + \sum\limits_{n\geq 1} \mathrm{ch}(\Gamma)u^n = \sum\limits_{n=\nu_1+2\nu_2+\ldots} \dfrac{\hbar^{\nu_1+\nu_2+\ldots}u^{\nu_1+2\nu_2+\ldots}}{1^{\nu_1}2^{\nu_2}\ldots\nu_1!\nu_2!\ldots}p_1^{\nu_1}p_2^{\nu_2}\ldots
\end{equation*}
Rewrite the last sum in more detail:
\begin{equation*}
\sum\limits_{n=\nu_1+2\nu_2+\ldots} \dfrac{\hbar^{\nu_1+\nu_2+\ldots}u^{\nu_1+2\nu_2+\ldots}}{1^{\nu_1}2^{\nu_2}\ldots\nu_1!\nu_2!\ldots}p_1^{\nu_1}p_2^{\nu_2}\ldots = \prod\limits_{n=1}^{\infty}\sum\limits_{k=0}^{\infty}\dfrac{\hbar^k u^{nk}}{n^kk!}p_n^k.
\end{equation*}
Note that this is the exponent
\begin{equation*}
\prod\limits_{n=1}^{\infty}\sum\limits_{k=0}^{\infty}\dfrac{\hbar^k u^{nk}}{n^kk!}p_n^k = \mathrm{exp}\,\hbar\,\left(up_1+\dfrac{u^2p_2}{2}+\dfrac{u^3p_3}{3}+\ldots\right),
\end{equation*}
and after rewriting the polynomials $p_i$ in the variables $x_i$ (see Eq.~(\ref{7})) we obtain
\begin{equation*}
\mathrm{exp}\,\hbar\,\left(up_1+\dfrac{u^2p_2}{2}+\dfrac{u^3p_3}{3}+\ldots\right) = \mathrm{exp}\,\hbar\,\sum\limits_{i=1}^\infty\left(ux_1+\dfrac{(ux_2)^2}{2}+\dfrac{(ux_3)^3}{3}+\ldots\right).
\end{equation*}
The last sum can be rewritten as
\begin{multline*}
\mathrm{exp}\,\hbar\,\sum\limits_{i=1}^\infty\left(ux_1+\dfrac{(ux_2)^2}{2}+\dfrac{(ux_3)^3}{3}+\ldots\right) = \mathrm{exp}\left(\,\hbar\,\sum\limits_{i=1}^\infty\mathrm{ln}(1-ux_i)\right) = \\
= \prod\limits_{i=1}^\infty (1-ux_i)^{-\hbar}.
\end{multline*}

Therefore, it remains to prove that
\begin{equation*}
\prod\limits_{i=1}^\infty (1-x_i)^{-\hbar} = \sum\limits_{\nu\vdash n} \prod\limits_{k\in\nu}\dfrac{\hbar+c(k)}{h(k)}\cdot s_\nu(x_1,x_2,\ldots).
\end{equation*}

The right hand side of the last equality is polynomial in~$\hbar$, whence it suffices to prove
the equality for natural values of $\hbar$. For $N\in\mb{N}$, we will prove the equality
\begin{equation}
\prod\limits_{i=1}^\infty (1-x_i)^{-N} = \sum\limits_{\nu\vdash n} \prod\limits_{k\in\nu}\dfrac{N+c(k)}{h(k)}\cdot s_\nu(x_1,x_2,\ldots).\label{2}
\end{equation}
We will use two statements from \cite{Mc}.
First,
\begin{equation}
\prod_{j=1}^{\infty}\prod_{i=1}^{\infty}(1-u_jx_i)^{-1} = \sum_\nu s_\nu(u_1,u_2,\ldots)s_\nu(x_1,x_2,\ldots).\label{3}
\end{equation}
The relationship between the left and the right hand sides of Eq.~(\ref{3}) follows
from statements about complete symmetric polynomials. The complete sym-\\metric polynomial~$h_r$ is defined by the formula
\begin{equation*}
h_r(x_1,x_2,\ldots) = \sum\limits_{|\nu|=r}\sum\limits_{\tau} x_1^{\nu_{\tau(1)}}x_2^{\nu_{\tau(2)}}\ldots x_{l(\nu)}^{\nu_{\tau(l(\nu))}},
\end{equation*}
where the second summation runs over all permutations~$\tau$ of the parts of the partition $\nu$.

The generating series $H(t)$ for the polynomials $h_r$ can be represented in the form
$H(t) = \sum\limits_{r\geq 0}h_rt^r = \prod\limits_{i\geq 1}(1-x_it)^{-1}.$ Note that we can represent
every factor on the right hand side of the last equality as the sum of an infinite geometric progression.
For details of the proof of Eq.~(\ref{3}), see \cite{Mc}, Sec.~1.4, p. 48.

Next, \cite{Mc}, Sec.~1.3, Example~4 yields
\begin{equation}
\prod\limits_{k\in\nu}\dfrac{N+c(k)}{h(k)} = s_\nu\underbrace{(1,\ldots,1)}_{N}.\label{4}
\end{equation}
Recall (Definition \ref{d1}) that  $s_\nu(x_1,x_2,\ldots)$ is the quotient of two determinants.
If we substitute $x_i = q^{i-1}$, then both determinants will become Vandermonde determinants and
\begin{equation}
s_\nu(1,q,q^2,\ldots) = q^{n(\nu)}\prod\limits_{x\in\nu}\dfrac{1-q^{n+c(x)}}{1-q^{h(x)}}, \label{6}
\end{equation}
where $n$ is the total amount of parts in the partition $\nu$ and $n(\nu) = \sum\limits_{i=1}^n (i-1)\nu_i$
(for details, see \cite{Mc}, Sec.~1.3, Example~1). Now, if we set $t=1$,
 then the expression $\dfrac{\prod\limits_{i=1}^k (1-t^i)}{(1-t)^k}$ is equal to $k!$. Hence, Eq.~(\ref{4}) follows from Eq.~(\ref{6}) after substituting $q=1.$

Equation~(\ref{2}) becomes obvious if we put $u_1=\ldots=u_N=1,\;u_{N+1}=\ldots = 0$ in Eq.~(\ref{3}) and use Eq.~(\ref{4}). The lemma is proved.
$\qed$

{\imp The eigenvalue $B_\nu$ of the operator $B$ on the eigenvector $\chi^\nu$ is
$$
B_\nu = \dfrac{\mathrm{dim}_\nu}{n!} \prod\limits_{k\in\nu}(1+c(k)\hbar).
$$}
\proof
Consider the function $\tilde\Gamma:\sigma_\nu\mapsto \hbar^{n-l(\nu)}$. Let us expand $\tilde\Gamma$ in the basis of characters using Lemma~\ref{l}:
\begin{multline*}
 \tilde\Gamma(\cdot)= \sum\limits_{\nu\vdash n} \hbar^n \prod\limits_{k\in\nu}\dfrac{\dfrac{1}{\hbar}+c(k)}{h(k)}\chi^\nu = \sum\limits_{\nu\vdash n} \dfrac{\hbar^n}{\hbar^n} \prod\limits_{k\in\nu}\dfrac{1+c(k)\hbar}{h(k)}\chi^\nu = \\
= \sum\limits_{\nu\vdash n} \dfrac{\mathrm{dim}_\nu}{n!}\prod\limits_{k\in\nu}(1+c(k)\hbar)\chi^\nu.
\end{multline*}
To deduce the last equality, we made use of the famous hook-length formula:
\begin{equation*}
\mathrm{dim}_\nu = \dfrac{n!}{\prod\limits_{k\in\nu} h(k)}.
\end{equation*}
The corollary is proved.
$\qed$
\\
\textit{Proof of Theorem~\ref{t1}.\hspace{2mm}}
Consider the series
\begin{equation*}
\sum_{n=0}^\infty\sum_{\nu\vdash n} \prod\limits_{k\in\nu}(1+c(k)\hbar)^m \dfrac{\mathrm{dim}_\nu}{n!} s\hbar_\nu(\hbar,p_1,p_2,\ldots) \hbar^{-2n}
\end{equation*}
in variables $\hbar$ and $p_1, p_2,\ldots$. Its coefficients are the eigenvalues of the operator $B^m$ on the eigenvectors $\chi^\nu$ written in the basis of shifted Schur functions. This power series enumerates both connected and disconnected coverings. Take the logarithm of this series to isolate connected coverings.

The Riemann--Hurwitz formula gives the genus of the covering surface:
\begin{equation*}
2g = 2 - 2n + \sum\limits_{P} (k(P)-1).
\end{equation*}
Here the summation is over all ramification points $P$ in $\mb{C}P^1$, and $k(P)$ is the ramification order at the ramification point $P$. Recall that the ramification order~$k(P)$ at the ramification point $P$ is
defined as follows.
Let $g_P$ be the monodromy permutation of the covering over the ramification point $P$ and let $l$ be the number of cycles in $g_P$. Then $k(P)=n-l$.
Notice that the contribution
to the power of~$\hbar$ from scaling of the Schur functions is equal to the sum $\sum\limits_P (k(P)-1)$,
  due to the Riemann--Hurwitz formula. We multiply by $\hbar^{-2n}$ and by $\hbar^2$ (see Eq.~(\ref{1}))
  in order to compensate the exponent of the variable $\hbar$ and the genus of the covering surface taken twice.
$\qed$

\vspace{4mm}
It would be interesting to give a geometric interpretation of these formulas. In particular, is it true that the
Bousquet-Melou--Schaeffer numbers are equal to the intersection numbers of certain natural
characteristic classes on suitable moduli spaces similar to the well-known ELSV formula \cite{ELSV}?
 Let us recall this formula enumerating coverings by a surface of given genus~$g$
 with simple ramification over all ramification points but one, where the ramification type
 is $\nu_1,\ldots,\nu_t$:
\begin{equation*}
h_{g,\{\nu_1,\ldots\nu_t\}} = \left(\prod_{i=1}^t \dfrac{\nu_i^{\nu_i}}{\nu_i!}\right)\int\limits_{\mathcal{\overline{M}}_{g,t}} \dfrac{1 - \l_1+ \ldots + (-1)^g\l_g}{(1-\nu_1\psi_1)\cdot\ldots\cdot(1-\nu_t\psi_t)}.
\end{equation*}

In particular, for the genus zero coverings it gives Hurwitz's formula \cite{Hu}:
\begin{equation*}
h_{0,\{\nu_1,\ldots\nu_t\}} = n^{c-3}\left(\prod_{i=1}^t \dfrac{\nu_i^{\nu_i}}{\nu_i!}\right).
\end{equation*}

The Bousquet-Melou--Schaeffer formula
\begin{equation*}
|\mathrm{Aut(\nu)}| b_{0,\nu,m} = m((m-1)n-1)_{t-3}\prod\limits_{i=1}^t \dfrac{(m\nu_i-1)_{\nu_i}}{\nu_i!}
\end{equation*}
is very similar to the Hurwitz formula for $h_{0,\{\nu_1,\ldots\nu_t\}}$, with powers replaced by the descending
factorials, but its geometric interpretation
is unknown.

\section{Acknowledgements}

The author is very grateful to S.~Lando for his assistance and permanent attention, to M.~Kazaryan,
who conjectured an explicit form of the eigenvalues for the operator acting by multiplication
by the sum of all permutations, and to G.~Olshanski for explaining his work about functions on symmetric groups.

\renewcommand{\refname}{References}

\end{document}